# Multivariate sequential analysis with linear boundaries

Robert Keener[1],*

*University of Michigan*

**Abstract:** Let $\{S_n = (X_n, W_n)\}_{n \geq 0}$ be a random walk with $X_n \in \mathbb{R}$ and $W_n \in \mathbb{R}^m$. Let $\tau = \tau_a = \inf\{n : X_n > a\}$. The main results presented are two term asymptotic expansions for the joint distribution of $S_\tau$ and $\tau$ and the marginal distribution of $h(S_\tau/a, \tau/a)$ in the limit $a \to \infty$. These results are used to study the distribution of $t$-statistics in sequential experiments with sample size $\tau$, and to remove bias from confidence intervals based on Anscombe's theorem.

## 1. Introduction

The protocol for most sequential experiments dictate that data are collected until a random walk $\{X_n\}_{n \geq 0}$ crosses a curved boundary. For design and analysis of these experiments there has been great interest in the joint distribution of the sample size $\tau$ and $X_\tau$, or related marginal distributions. Using nonlinear renewal theory or other techniques, it is often possible to approximate distributions for experiments with curved boundaries by distributions when the boundary is linear. The results presented here all concern the linear case with sample size

$$\tau = \tau_a = \inf\{n : X_n > a\}.$$

As $a \to \infty$, the sample size $\tau$ is asymptotically normal (after suitable centering and scaling); the excess over the boundary, $S_\tau - a$, has a limiting distribution; and these variables are asymptotically independent. For this result, its extension to curved boundaries, and applications to sequential testing and estimation, see [24] or [17].

In applications, $a$ may not be very large, and to obtain an adequate approximation it is often necessary to refine the basic limit theory. In univariate situations, asymptotic expansions (with curved boundaries) are given in [9, 11, 18, 19, 25, 26, 29, 30].

In this paper, the primary concern will be with distribution theory for multivariate experiments in which several variables are recorded for each experimental unit. In practice, the stopping time will often be driven by the values for some primary variable, but parameters associated with distributions of the other variables will also be of interest. In a clinical setting, the primary variable might be some measure of the efficacy of treatment, with other variables related to potential side effects. Since the variables not directly related to stopping may be correlated with the primary stopping variable, a statistical analysis that ignores effects of optional stopping may be biased, even if the variables involved are not directly related to the stopping rule. Whitehead [22] notes this possibility and provides adjustments to reduce the bias of maximum likelihood estimators in these situations.

*Research supported in part by NSF Grant DMS-85-04708

[1]Department of Statistics, University of Michigan, 1085 South University, Ann Arbor, MI 48109, USA, e-mail: keener@umich.edu

Most of the expansions cited above for univariate data are based on arguments that do not extend easily to higher dimensions. A notable exception would be an approach introduced by Woodroofe [26]. This approach is based on the likelihood function in a parametric context, making it particularly appropriate for sequential or adaptive experiments since likelihood functions are not effected by optional stopping. With this approach, coverage probabilities for approximate confidence intervals are viewed as functions of the unknown parameter. Expansions for these coverage probabilities do not hold in the conventional pointwise sense. Instead, the expansions hold if the functions are smoothed by integration against some weight function. Woodroofe calls these results "very weak expansions." Very weak expansions have now been used in a variety of situations to set confidence intervals for secondary parameters—see [6, 13, 20, 21, 22, 23, 27, 28].

One final approach to interval estimation deserves mention. In one dimension, Siegmund [15] suggests ordering points in the stopping region and inverting a family of tests. Probability calculations in this approach are based on large deviation approximations, which gives reasonable accuracy in many standard designs—especially repeated significance tests. In higher dimensions, it is less evident how points should be ordered in the stopping region. However, Siegmund in [16] successfully extends his method to interval estimation following a repeated significance test of several normal means. Unfortunately, his argument uses the spherical symmetry of the normal distribution and may not be widely applicable.

In this paper, expansions are derived from related expansions for multivariate renewal measures given in [12]. In contrast with very weak expansions, a parametric model and likelihood function are not necessary, and the expansions hold in the conventional sense—no smoothing is necessary. Unfortunately, the algebra necessary with this approach seems more difficult. Expansions are only obtained with a leading correction term, and, as mentioned above, the stopping time $\tau$ must have a linear boundary. The main results give expansions for joint distributions of partial sums and the stopping time. Using these results, refined approximations for the distributions of $t$-statistics are derived and used to set confidence intervals for the mean of a secondary variables. The refined distribution theory reduces tail probability bias present in confidence intervals set using the normal approximation based on Ancsombe's theorem.

Extensions of the results in this paper to curved boundaries seem challenging but may be possible. In the univariate case, considered in [11], under conditions sufficient for local limit theorems the curvature of the boundary only has an intrinsic effect on the second correction term—the chance of crossing a curved boundary on step $n$ agrees with the chance of crossing an appropriate linear boundary at step $n$ to $o(1/a)$. But curvature does have some effect on the first correction term since the appropriate linear boundary must vary with $n$. Using this, it may be possible to guess how the results here could be modified for curved boundaries, although the best approach for a rigorous argument is not clear.

Let $\{(X_n, W_n)\}_{n \geq 0}$ be a random walk with $X_n \in \mathbb{R}$ and $W_n \in \mathbb{R}^m$. Theorem 1.1 provides an expansions for the joint distribution of $X_\tau$ and $W_\tau$ and the joint distribution of $X_\tau$, $W_\tau$ and $\tau$. To state this theorem, a fair amount of notation is necessary. Let $\phi$ denote the normal density with mean zero and identity covariance. For any function $f$, the *oscillation function* $\omega_f$ is given by

$$\omega_f(x; \epsilon) = \sup\{|f(x) - f(y)| : |x - y| \leq \epsilon\}.$$

The random vector $(X, W) = (X_1, W_1)$ satisfies Cramér's condition if

$$\limsup_{(\xi_1, \xi_2) \to \infty} |Ee^{i(\xi_1 X + \xi_2 \cdot W)}| < 1.$$

For joint expansions with $\tau$, let $W_n^* = (W_n, n)$ and take $W^* = W_1^*$. Then $W_\tau^* = (W_\tau, \tau)$. Let $\nu = EX$, $\gamma = EW/\nu$, $\gamma^* = EW^*/\nu$, $\Sigma = \text{Cov}(W - \gamma X)$, $\Sigma_* = \text{Cov}(W^* - \gamma^* X)$, $Z_n = \Sigma^{-1/2}(W_n - \gamma X_n)$, $Z_n^* = \Sigma_*^{-1/2}(W_n^* - \gamma^* X_n)$, $Z = Z_1$ and $Z^* = Z_1^*$. The first ladder time is $T = \tau_0$. Let $\tilde{X} = X_T$, $\tilde{W} = W_T$, $\tilde{W}^* = W_T^*$, $\tilde{Z} = Z_T/\sqrt{ET}$ and $\tilde{Z}^* = Z_T^*/\sqrt{ET}$. Define

$$\rho_0(x) = \frac{P(\tilde{X} \geq x)}{\nu ET},$$

$$\rho_1(x) = \frac{E[\tilde{Z}; \tilde{X} \geq x]}{\nu\sqrt{ET}}$$

and

$$\rho_1^*(x) = \frac{E[\tilde{Z}^*; \tilde{X} \geq x]}{\nu\sqrt{ET}}$$

for $x > 0$, with $\rho_0$, $\rho_1$ and $\rho_1^*$ identically zero on $(-\infty, 0]$. The densities for the approximate distributions in Theorem 1.1 are given by

$$\frac{d\hat{Q}}{d\lambda}(x, w) = \frac{\phi(q)}{\sqrt{|\Sigma|}(a/\nu)^{m/2}} \Big\{ \rho_0(x) + \sqrt{\nu/a} \big[ \mathbf{H}(q)\rho_0(x) + q \cdot \rho_1(x) \big] \Big\}$$

and

$$\frac{d\hat{Q}^*}{d\lambda^*}(x, w^*) = \frac{\phi(q^*)}{\sqrt{|\Sigma_*|}(a/\nu)^{(m+1)/2}} \Big\{ \rho_0(x) + \sqrt{\nu/a} \big[ \mathbf{H}^*(q^*)\rho_0(x) + q^* \cdot \rho_1^*(x) \big] \Big\},$$

where $\lambda$ is Lebesgue measure on $\mathbb{R}^{m+1}$, $\lambda^*$ is the product of $\lambda$ with counting measure on $\mathbb{Z}$,

$$q = \Sigma^{-1/2}(w - \gamma x - \gamma a)\sqrt{\nu/a},$$

$$q^* = \Sigma_*^{-1/2}(w^* - \gamma^* x - \gamma^* a)\sqrt{\nu/a},$$

$$\mathbf{H}(q) = \frac{1}{6}E(q \cdot Z)^3 - \frac{1}{2}EZ^2 q \cdot Z$$
$$+ \frac{(m + 2 - |q|^2)EXq \cdot Z}{2\nu} - \frac{E\tilde{X}q \cdot \tilde{Z}}{\nu\sqrt{ET}}$$

and

$$\mathbf{H}^*(q^*) = \frac{1}{6}E(q^* \cdot Z^*)^3 - \frac{1}{2}EZ^{*2}q^* \cdot Z^*$$
$$+ \frac{(m + 3 - |q^*|^2)EXq^* \cdot Z^*}{2\nu} - \frac{E\tilde{X}q^* \cdot \tilde{Z}^*}{\nu\sqrt{ET}}.$$

**Theorem 1.1.** *Suppose $(X, W)$ satisfies Cramér's condition, $\nu > 0$, $E|X|^{(3+\delta)/2} < \infty$ and $E|Z|^{3+\delta} < \infty$, where $\delta \in (0, 1)$. Then for some $\eta > 0$,*

$$Ef(X_\tau - a, W_\tau) = \int f \, d\hat{Q} + O(1) \int \omega_f(\cdot; e^{-\eta a}) \, d\hat{Q}$$
$$+ o\{a^{(-1-\delta)/2}(\log a)^{m/2}\}$$

as $a \to \infty$, uniformly for measurable nonnegative $f$ bounded above by 1. If the moment condition for $X$ is strengthened to $E|X|^{3+\delta} < \infty$, then for some $\eta > 0$,

$$Ef(X_\tau - a, W_\tau^*) = \int f \, d\hat{Q}^* + O(1) \int \omega_f(\cdot; e^{-\eta a}) \, d\hat{Q}^*$$
$$+ o\{a^{(-1-\delta)/2}(\log a)^{(m+1)/2}\}$$

as $a \to \infty$, uniformly for measurable nonnegative $f$ bounded above by 1.

With $f$ an indicator function, this theorem provides approximations for probabilities. Unfortunately, access to the wealth of information available in principle from this result may be rather difficult. The next result uses Theorem 1.1 to approximate the marginal distribution of $h(X_\tau/a, W_\tau/a, \tau/a)$ for smooth functions $h$. Special cases of interest include linear functions, averages, normalized partial sums and $t$-statistics.

To state this result, let $Y_n$ denote the first coordinate of $W_n$, so $W_n = (Y_n, V_n)$ where $V_n \in \mathbb{R}^{m-1}$. Also, let $V_n^* = (V_n, n)$, so $W_n^* = (Y_n, V_n^*)$. Finally, let $S_n^* = (S_n, n) = (X_n, Y_n, V_n, n)$. By the strong law of large numbers, $X_\tau/a \to 1$, $Y_\tau/a \to \gamma_1$ and $V_\tau^*/a \to \gamma_2$ as $a \to \infty$, where $\gamma_1 = EY/\nu$ and $\gamma_2 = EV^*/\nu$. The expansion naturally involves Taylor expansion of $h$ about $s_0 = (1, \gamma_1, \gamma_2)$. Let $\nabla_{v^*}$ and $\nabla^2_{v^*}$ denote the gradient vector and Hessian matrix with respect to $v^*$. The regularity assumptions for $h$ are as follows:

1. The mixed third derivatives of $h$ are continuous in some neighborhood $N_0$ of $s_0$ (this assumption could be relaxed slightly).
2. At $s_0$, $h = \partial h/\partial x = 0$, $\nabla_{v^*} h = 0$ and $\partial h/\partial y = 1$

Since invertible affine transformations of $(X, W)$ preserve the moment conditions and Cramér's condition, this last assumption is less restrictive than it may appear. The only cases not covered are those where the first order Taylor expansion of $h(S_\tau^*/a)$ is linearly independent of $W_\tau$, and hence depends only on $X_\tau$ and $\tau$. Define

$$h_0 = \frac{1}{2} \frac{\partial^2}{\partial y^2} h(s_0),$$

$$A = \frac{1}{2} \nabla^2_{v^*} h(s_0),$$

and let $h_1$ be the linear function given by

$$h_1(v^*) = v^* \cdot \nabla_{v^*} \frac{\partial}{\partial y} h(s_0).$$

The quadratic Taylor approximation for $h$ about $s_0$ is

$$(y - \gamma_1) + h_0(y - \gamma_1)^2 + (y - \gamma_1)h_1(v^* - \gamma_2)$$
$$+ (v^* - \gamma_2) \cdot A(v^* - \gamma_2) + h_3(x - 1)^2 + (x - 1)h_4(y - \gamma_1, v^* - \gamma_2),$$

where $h_3$ is a constant and $h_4$ is a bi-linear function—their specification is not important. Partition $\Sigma_*$ as

$$\Sigma_* = \begin{pmatrix} \Sigma_{11} & \Sigma_{12} \\ \Sigma_{21} & \Sigma_{22} \end{pmatrix}$$

where $\Sigma_{22} = \mathrm{Cov}(V^* - \gamma_2 X)$ is $m$ by $m$ and $\Sigma_{11} = \mathrm{Var}(Y - \gamma_1 X) = \sigma^2$. The approximate distribution function in Theorem 1.2 is

$$F_a(c) = \Phi(\hat{c}) + \phi(\hat{c})\sqrt{\frac{\nu}{a}}\bigg\{\bigg[-\frac{1}{6}E(Y-\gamma_1 X)^3 + \frac{EX(Y-\gamma_1 X)}{2\nu}\bigg](\hat{c}^2 - 1) - \frac{\gamma_1 E\tilde{X}^2}{2\nu\sigma ET}$$
$$- \frac{\sigma\hat{c}^2 h_0}{\nu} - \frac{\hat{c}^2 h_1(\Sigma_{21})}{\nu\sigma} - \frac{\hat{c}^2 - 1}{\nu\sigma^3}\Sigma_{12} A \Sigma_{21} - \frac{\mathrm{tr}[A\Sigma_{22}]}{\nu\sigma}\bigg\},$$

where $\hat{c} = c\sqrt{\nu}/\sigma$.

**Theorem 1.2.** *Suppose $0 < \delta < 1$, $E|X|^{3+\delta} < \infty$, $E|W|^{3+\delta} < \infty$ and $(X, W)$ satisfies Cramér's condition. Also, if $\delta > \sqrt{2} - 1$, assume $E|X|^{2/(1-\delta)} < \infty$. Let*

$$\Xi = \sqrt{a}h(X_\tau/a, W_\tau/a, \tau/a),$$

*where $h$ satisfies the conditions stated above. Then*

$$P(\Xi \leq c) = F_a(c) + o\big\{a^{(-1-\delta)/2}(\log a)^{(m+1)/2}\big\}$$

*as $a \to \infty$, uniformly for $c \in \mathbb{R}$.*

The rest of the paper is organized as follows. In the next section, Theorem 1.2 is specialized to $t$-statistics and used to set confidence intervals for $EY$. Results are reported for a simulation study comparing the coverage probabilities of these confidence intervals with confidence intervals set using Anscombe's theorem. Section 3 contains an approximation for the joint distribution of $X_\tau$ and $W_\tau$ in the positive case where $X > 0$. In Section 4, ladder variables are introduced and used to prove Theorem 1.1 from results in Section 3. In Section 5, various marginal distributions are approximated and Theorem 1.2 is proved.

## 2. Confidence intervals and $t$-statistics

The main concern of this section is setting confidence intervals for $\mu = EY$ after a sequential experiment with sample size $\tau$. Let $e_n = Y_n - Y_{n-1}$, so $Y_n = \sum_{i=1}^n e_i$. The $t$-statistic for $\mu$ is

$$T = \frac{\overline{Y}_\tau - \mu}{\hat{\sigma}/\sqrt{\tau}},$$

where $\overline{Y}_\tau = Y_\tau/\tau$ and

$$\hat{\sigma}^2 = \frac{1}{\tau - 1}\sum_{i=1}^\tau (e_i - \overline{Y}_\tau)^2.$$

By Anscombe's theorem [1],

(1) $$T \Rightarrow N(0, 1)$$

as $a \to \infty$. Using this, the coverage probability of the confidence interval

$$(\mathrm{LCL}_0, \mathrm{UCL}_0) = (\overline{Y}_\tau \pm z_\alpha \hat{\sigma}/\sqrt{\tau})$$

approaches $1 - 2\alpha$ as $a \to \infty$, where $z_\alpha = \Phi^{-1}(1 - \alpha)$.

The normal approximation (1) can be improved using Theorem 1.2. It is most convenient to work initially with the modified statistic

$$T_0 = \frac{\overline{Y}_\tau - \mu}{\hat{\sigma}_0/\sqrt{\tau}},$$

where $\hat{\sigma}_0^2 = (\tau - 1)\hat{\sigma}^2/\tau$. If $V_n = \sum_{i=1}^n e_i^2$, then

$$\hat{\sigma}_0^2 = \left\{ \frac{V_\tau}{\tau} - \frac{Y_\tau^2}{\tau^2} \right\}$$

and

$$T_0 = \sqrt{a} h(Y_\tau/a, V_\tau/a, \tau/a),$$

where

$$h(y, v, t) = \frac{y - t\mu}{\sqrt{v - y^2/t}}.$$

**Corollary 2.1.** *If $(X, Y, Y^2)$ satisfies Cramér's condition, $\nu > 0$, $E|X|^{(3+\delta)/2} < \infty$, $E|X|^{2/(1-\delta)} < \infty$ and $E|Y|^{6+2\delta} < \infty$, then*

$$P(T_0 \leq c) \approx \Phi(c) + \phi(c)\sqrt{\frac{\nu}{a}} \left\{ \frac{\mu_3}{6\sigma^3}(1 + 2c^2) - \frac{\Sigma_{XY}}{2\nu\sigma} \right\}$$

*uniformly in $c$ as $a \to \infty$, where $\sigma^2 = \mathrm{Var}(Y)$, $\mu_3 = E(Y - \mu)^3$, $\Sigma_{XY} = \mathrm{Cov}(X, Y)$. Here and later in this section, "$\approx$" will mean that the two quantities differ by $o\{a^{(-1-\delta)/2}(\log a)^{3/2}\}$.*

*Proof.* After recentering $Y$ and rescaling $X$, $Y$ and $a$, there is no harm assuming $\mu = 0$, $\sigma = 1$ and $\nu = 1$. Then $\hat{c} = c$, $\gamma_1 = 0$ and $\gamma_2 = (1, 1)$. Differentiation gives $h_0 = 0$, $h_2 = 0$ and $h_1(v, t) = -v/2$, so

$$F_a(c) = \Phi(c) + \frac{\phi(c)}{\sqrt{a}} \left\{ \left[ -\frac{1}{6} EY^3 + \frac{1}{2} EXY \right](c^2 - 1) + \frac{c^2}{2} E(Y^3 - XY) \right\}$$

and the corollary follows from Theorem 1.2. □

The next lemma is a generalization of Slutsky's theorem to asymptotic expansions. Using it, the approximation in this corollary is also an approximation for the distribution of $T$. By definition, say that $\zeta_a \to 0$ in probability at rate $o(\delta_a)$ if for every $\epsilon > 0$,

$$P(|\zeta_a| > \epsilon) = o(\delta_a)$$

as $a \to \infty$.

**Lemma 2.2.** *Suppose the distributions of a family of variables $\{\xi_a\}$ have an asymptotic expansion of the form*

$$P(\xi_a \leq c) = G_a(c) + o(\epsilon_a)$$

*as $a \to \infty$, uniformly in $c$, where the $G_a$ satisfy*

$$\limsup_{a \to \infty} \sup_x \sup_{y > 0} \{G_a(x + y) - G_a(x)\}/y < \infty.$$

*Then if $\zeta_a \to 0$ in probability at rate $o(\epsilon_a)$,*

$$P(\xi_a + \epsilon_a \zeta_a \leq c) = G_a(c) + o(\epsilon_a)$$

*as $a \to \infty$, uniformly in $c$.*

*Proof.* Let
$$K = \sup_{a > a_0} \sup_x \sup_{y > 0} \{G_a(x+y) - G_a(x)\}/y.$$

Then $K < \infty$ if $a_0$ is large enough. For any $\epsilon > 0$, for $a > a_0$,

$$\begin{aligned} P(\xi_a + \epsilon_a \zeta_a \leq c) &\leq P(\xi_a \leq c + \epsilon \epsilon_a) + P(\zeta_a > \epsilon) \\ &= G_a(c + \epsilon \epsilon_a) + o(\epsilon_a) \\ &\leq G_a(c) + K \epsilon \epsilon_a + o(\epsilon_a) \end{aligned}$$

as $a \to \infty$, uniformly in $c$. Since $\epsilon$ is arbitrary,

$$P(\xi_a + \epsilon_a \zeta_a \leq c) \leq G_a(c) + o(\epsilon_a)$$

as $a \to \infty$, uniformly in $c$. The reverse inequality follows in a similar fashion. □

Using Lemma 2.2, it is easy to check that $T$ and $T_0$ have the same asymptotic expansion (to this order). For setting confidence intervals, it is convenient (following Hall [10]) to write the expansion in the form

$$P\left(T \leq c - \sqrt{\frac{\nu}{a}} \left[\frac{\mu_3}{6\sigma^3}(1 + 2c^2) - \frac{\Sigma_{XY}}{2\nu\sigma}\right]\right) \approx \Phi(c)$$

as $a \to \infty$. This expansion may not hold uniformly in $c$ as $a \to \infty$, but will hold uniformly for $c$ in any compact set. It is now natural to replace the parameters $\nu$, $\sigma$, $\mu_3$ and $\Sigma_{XY}$ by estimates $\hat{\nu}$, $\hat{\sigma}$, $\hat{\mu}_3$ and $\hat{\Sigma}_{XY}$. By the same calculations used to prove Lemma 2.2, if

$$a^{\delta/2}(\log a)^{3/2}\{|\hat{\nu} - \nu| + |\hat{\sigma} - \sigma| + |\hat{\mu}_3 - \mu_3| + |\hat{\Sigma}_{XY} - \Sigma_{XY}|\} \to 0$$

in probability at rate $o\{a^{(-1-\delta)/2}(\log a)^{3/2}\}$, then

$$P\left(T \leq c - \sqrt{\frac{\hat{\nu}}{a}} \left[\frac{\hat{\mu}_3}{6\hat{\sigma}^3}(1 + 2c^2) - \frac{\hat{\Sigma}_{XY}}{2\hat{\nu}\hat{\sigma}}\right]\right) \approx \Phi(c)$$

as $a \to \infty$, uniformly for $c$ in any compact set. With $c = \pm z_\alpha$, the coverage probability of the confidence interval

$$(\text{LCL}_1, \text{UCL}_1) = \overline{Y}_\tau + \frac{\hat{\sigma}}{\sqrt{\tau}}\sqrt{\frac{\hat{\nu}}{a}}\left\{\frac{\hat{\mu}_3}{6\hat{\sigma}^3}(1 + 2z_\alpha^2) - \frac{\hat{\Sigma}_{XY}}{2\hat{\nu}\hat{\sigma}}\right\} \pm z_\alpha \frac{\hat{\sigma}}{\sqrt{\tau}}$$

is $\approx 1 - 2\alpha$ as $a \to \infty$.

In parametric estimation following a sequential test, Whitehead [22] suggests centering confidence intervals around the maximum likelihood estimator less an estimate of its bias. In the normal case where $\overline{Y}_\tau$ is the maximum likelihood estimator for $\mu$, the interval suggested here is different: the center to the relevant order is $\overline{Y}_\tau$ less *half* its bias.

Corollary 2.1 can be used to study the performance of the interval $(\text{LCL}_0, \text{UCL}_0)$. By the corollary,

$$P(\mu \geq \text{UCL}_0) \approx \alpha + \phi(z_\alpha)\sqrt{\frac{\nu}{a}}\left\{\frac{\mu_3}{6\sigma^3}(1 + 2z_\alpha^2) - \frac{\Sigma_{XY}}{2\nu\sigma}\right\}$$

and
$$P(\mu \leq \text{LCL}_0) \approx \alpha - \phi(z_\alpha)\sqrt{\frac{\nu}{a}}\left\{\frac{\mu_3}{6\sigma^3}(1+2z_\alpha^2) - \frac{\Sigma_{XY}}{2\nu\sigma}\right\},$$

so
$$P(\text{LCL}_0 < \mu < \text{UCL}_0) \approx 1 - 2\alpha.$$

To this order of analysis, the overall coverage probability equals the desired value, $1 - 2\alpha$, but this confidence interval is biased—the two error probabilities, $P(\mu \geq \text{UCL}_0)$ and $P(\mu \leq \text{LCL}_0)$ do not agree. The modified interval $(\text{LCL}_1, \text{UCL}_1)$ has no bias to this order of analysis:

$$P(\mu \geq \text{UCL}_1) \approx \alpha \approx P(\mu \leq \text{LCL}_1).$$

Although the modified confidence interval has less bias and comparable coverage as $a \to \infty$, there is the usual concern that these properties may fail for moderate values of $a$. Extra variation is introduced when the parameters $\nu$, $\sigma$, $\mu_3$ and $\Sigma_{XY}$ are estimated, and this may degrade the overall coverage of the interval. The performance of several confidence intervals are studied in the simulation study reported in Table 1. In this study, $(X, Y)$ follow a bivariate normal distribution with $\text{Var}(X) = \text{Var}(Y) = 1$, so $\Sigma_{XY}$ equals the correlation coefficient $\rho$. Two versions of the modified confidence interval are considered. For both, $\nu$ and $\Sigma_{XY}$ are estimated as

$$\hat{\nu} = \overline{X}_\tau = X_\tau/\tau$$

and
$$\hat{\Sigma}_{XY} = \frac{1}{\tau}\sum_{i=1}^{\tau}(X_i - X_{i-1})e_i - \overline{X}_\tau\overline{Y}_\tau.$$

For the interval $(\text{LCL}_1^{(1)}, \text{UCL}_1^{(1)})$, $\hat{\mu}_3 = 0$. This would be appropriate in practice if the researcher were fairly certain the centered marginal distribution of $Y$ was symmetric. For the interval $(\text{LCL}_1^{(2)}, \text{UCL}_1^{(2)})$, $\mu_3$ is estimated by

$$\hat{\mu}_3 = \frac{1}{\tau}\sum_{i=1}^{\tau}(e_i - \overline{Y}_\tau)^3.$$

In the simulation, $\alpha$ was fixed at 5%. There were 10,000 replications.

Examining the simulation results a few tentative conclusions seem in order. There is almost no difference between the performance of the two versions of the modified confidence intervals. The original confidence interval has better coverage, but only slightly—less than one percent. The original confidence interval can be quite biased with more that a 2 to 1 ratio between the error probabilities in some cases. The modified intervals were always less biased. Finally, it is worth noting that in all cases the overall coverage is a bit less less than 90%. Higher order corrections would be useful.

## 3. Positive case

The notation in this section parallels the notation in the introduction, but is somewhat divergent. Let $\{(X_n, W_n)\}_{n\geq 0}$ be a random walk with $X = X_1 > 0$, and $W = W_1 \in \mathbb{R}^m$. The stopping time is still $\tau = \tau_a = \inf\{n : X_n > a\}$. The main result of this section is an asymptotic expansion for the joint distribution of $X_\tau$ and $W_\tau$. When this result is used proving Theorem 1.1, $(X, W)$ will be either $(\tilde{X}, \tilde{W})$ or $(\tilde{X}, \tilde{W}^*)$. Since the last coordinate of $\tilde{W}^*$ is $T$, two different smoothness conditions will be used. In the "continuous" case, $(X, W)$ will satisfy Cramér's condition. In the "mixed" case, with $W = (Y, T)$ say, $T$ will be *arithmetic* on $\mathbb{Z}$, i.e., $P(T \in \mathbb{Z}) = 1$ but $P(T \in B) < 1$ for $B$ any proper subgroup of $\mathbb{Z}$, and $(X, Y)$ will be *strongly nonlattice* with $T$, i.e.,

$$\liminf_{\xi_1^2 + \xi_2^2 \to \infty} \inf_{-\pi < \xi_3 < \pi} \left| 1 - Ee^{i\xi_1 X + i\xi_2 \cdot Y + i\xi_3 T} \right| > 0.$$

TABLE 1
*Coverage Probabilities and Error Probabilities*

|  | $a = 10$ | | | | $a = 25$ | | | |
|---|---|---|---|---|---|---|---|---|
|  | $\nu = 0.5$ | | $\nu = 0.25$ | | $\nu = 0.5$ | | $\nu = 0.25$ | |
|  | $\rho = 0.4$ | $\rho = 0.8$ | $\rho = 0.4$ | $\rho = 0.8$ | $\rho = 0.4$ | $\rho = 0.8$ | $\rho = 0.4$ | $\rho = 0.8$ |
| $(\text{LCL}_0, \text{UCL}_0)$ | 0.884 | 0.881 | 0.892 | 0.894 | 0.890 | 0.887 | 0.890 | 0.893 |
| $P(\mu \geq \text{UCL}_0)$ | 0.047 | 0.035 | 0.045 | 0.030 | 0.051 | 0.047 | 0.050 | 0.042 |
| $P(\mu \leq \text{LCL}_0)$ | 0.069 | 0.083 | 0.064 | 0.075 | 0.059 | 0.066 | 0.061 | 0.065 |
| $\text{LCL}_1^{(1)}, \text{UCL}_1^{(1)})$ | 0.880 | 0.876 | 0.892 | 0.879 | 0.890 | 0.884 | 0.883 | 0.885 |
| $P(\mu \geq \text{UCL}_1^{(1)})$ | 0.057 | 0.055 | 0.057 | 0.062 | 0.057 | 0.061 | 0.064 | 0.062 |
| $P(\mu \leq \text{LCL}_1^{(1)})$ | 0.063 | 0.069 | 0.051 | 0.060 | 0.054 | 0.055 | 0.054 | 0.053 |
| $(\text{LCL}_1^{(2)}, \text{UCL}_1^{(2)})$ | 0.874 | 0.873 | 0.892 | 0.876 | 0.888 | 0.885 | 0.882 | 0.885 |
| $P(\mu \geq \text{UCL}_1^{(2)})$ | 0.058 | 0.056 | 0.059 | 0.063 | 0.058 | 0.060 | 0.063 | 0.061 |
| $P(\mu \leq \text{LCL}_1^{(2)})$ | 0.068 | 0.071 | 0.050 | 0.060 | 0.054 | 0.056 | 0.055 | 0.054 |

In the continuous case, $\lambda$ will denote Lebesgue measure on $\mathbb{R}^{m+1}$, and in the mixed case $\lambda$ will denote the product of Lebesgue measure on $\mathbb{R}^m$ with counting measure on $\mathbb{Z}$. Also, define $\tilde{\lambda}$ by

$$\tilde{\lambda}(B) = \lambda\bigl(\{(x, w) : (x, \Sigma^{-1/2}(w - \gamma x)) \in B\}\bigr)$$

for Borel sets $B$.

Let $\nu = EX$, $\gamma = EW/\nu$ and $\Sigma = \text{Cov}(W - \gamma X)$. Assume $0 < \Sigma < \infty$ and take $Z_n = \Sigma^{-1/2}(W_n - \gamma X_n)$ and $Z = Z_1$. Let $(X_n, W_n) \sim F_n$ and $(X_n, Z_n) \sim \tilde{F}_n$ with $F = F_1$ and $\tilde{F} = \tilde{F}_1$. Then $R = \sum_{n \geq 0} F_n$ and $\tilde{R} = \sum_{n \geq 0} \tilde{F}_n$ are the renewal measures for the random walks $\{(X_n, W_n)\}_{n \geq 0}$ and $\{(X_n, Z_n)\}_{n \geq 0}$ respectively.

**Theorem 3.1.** *Suppose $0 < \delta < 1$, $E|X|^{(3+\delta)/2} < \infty$, $E|Z|^{3+\delta} < \infty$ and $\Sigma$ is positive definite. Also, assume that either $(X, W)$ satisfies Cramér's condition, or $W = (Y, T)$ with $T$ arithmetic on $\mathbb{Z}$ and $(X, Y)$ strongly nonlattice with $T$. Then for some $\eta > 0$,*

$$Ef(X_\tau - a, W_\tau) = \int f\, d\hat{Q} + O(1) \int \omega_f(\cdot\,; e^{-\eta a})\, d\hat{Q} + o\bigl\{a^{(-1-\delta)/2}(\log a)^{m/2}\bigr\}$$

*as $a \to \infty$, uniformly for measurable nonnegative functions $f$ bounded above by 1, where for $x \geq 0$,*

$$\frac{d\hat{Q}}{d\tilde{\lambda}}(x, w) = \frac{\phi(q)}{\sqrt{|\Sigma|}(a/\nu)^{m/2}} \Bigl\{\rho_0(x) + \sqrt{\nu/a}\bigl[\mathbf{H}(q)\rho_0(x) - q \cdot \rho_1(x)\bigr]\Bigr\},$$

$$q = \Sigma^{-1/2}(w - \gamma x - \gamma a)\sqrt{\nu/a},$$

$$\mathbf{H}(q) = \frac{1}{6} E(q \cdot Z)^3 - \frac{1}{2} E Z^2 q \cdot Z + \frac{m - |q|^2}{2\nu} E X q \cdot Z,$$

$$\rho_0(x) = P(X \geq x)/\nu,$$

and

$$\rho_1(x) = E[Z; X \geq x]/\nu.$$

The expectation of interest here can be expressed as an integral against the renewal measure $\tilde{R}$. Since $X > 0$, the event $\tau = n+1$ is the same as $X_n \leq a$ and $X_{n+1} > a$. Hence

$$\begin{aligned}
Ef(X_\tau - a, Z_\tau) &= \sum_{n=0}^{\infty} E\big[f(X_{n+1} - a, Z_{n+1}); \tau = n+1\big] \\
&= \sum_{n=0}^{\infty} \iint_{a-x_0 < x \leq a} f(x + x_0 - a, z + z_0) \, d\tilde{F}_n(x,z) \, d\tilde{F}(x_0, z_0) \\
&= \iint_{a-x_0 < x \leq a} f(x + x_0 - a, z + z_0) \, d\tilde{R}(x,z) \, d\tilde{F}(x_0, z_0).
\end{aligned} \quad (2)$$

Expansions for multivariate renewal measures are studied in [12]. To proceed we will use the following result which follows almost immediately from Theorem 3 of [12].

**Theorem 3.2.** *Under the assumptions of Theorem 3.1, for some $\eta > 0$,*

$$\int_{a < x \leq a + \Delta} f(x, z) \, d\tilde{R}(x, z) = \int_{a < x \leq a + \Delta} \big[f(x, z) + O(1)\omega_f(x, z; e^{-\eta a})\big] \, d\hat{\tilde{R}}(x, z) + o\big\{(1+\Delta) a^{(-1-\delta)/2} (\log a)^{m/2}\big\}$$

*as $a \to \infty$, uniformly for $\Delta > 0$ and nonnegative measurable $f$ bounded by 1, where*

$$\frac{d\hat{\tilde{R}}}{d\tilde{\lambda}}(x, z) = \hat{\tilde{r}}(x, z) = \frac{\phi(z\sqrt{\nu/x})}{\nu\sqrt{|\Sigma|}(x/\nu)^{m/2}} \big\{1 + \sqrt{\nu/x}\mathbf{H}(z\sqrt{\nu/x})\big\}.$$

*Proof of Theorem 3.1.* The integral in (2) over $x \geq a/2$ or $|z| \geq \sqrt{a}$ as $a \to \infty$ is $o\{a^{(-1-\delta)/2}\}$. To see this, note first that

$$\int_{a-x_0 < x \leq a} d\tilde{R}(x, z) = \tilde{R}\big((a - x_0, a] \times \mathbb{R}^m\big)$$

which is bounded by a multiple of $1 + x_0$, uniformly in $a$ and $x_0 \geq 0$ by the renewal theorem in one dimension. Then

$$\begin{aligned}
\int_{x_0 \geq a/2} (1 + x_0) d\tilde{F}(x_0, z_0) &= E[1 + X; X \geq a/2] \\
&\leq E\big[(1 + X)^{(3+\delta)/2}; X \geq a/2\big](1 + a/2)^{(-1-\delta)/2} \\
&= o\{a^{(-1-\delta)/2}\}
\end{aligned} \quad (3)$$

and

$$\int_{|z_0| \geq \sqrt{a}} (1+x_0) d\tilde{F}(x_0, z_0) = E[1+X; |Z| \geq \sqrt{a}]$$

(4)
$$\leq E[(1+X)|Z|^{1+\delta}; |Z| \geq \sqrt{a}] a^{(-1-\delta)/2}$$
$$= o\{a^{(-1-\delta)/2}\}$$

as $a \to \infty$. We can now use the approximation of Theorem 3.2 in (2). This gives

$$Ef(X_\tau - a, Z_\tau)$$
$$= \iint_{\substack{a-x_0 < x \leq a \\ x_0 < a/2 \\ |z_0| < \sqrt{a}}} [f(x+x_0-a, z+z_0) + O(1)\omega_f(x+x_0-a, z+z_0; e^{-\eta a})]$$
$$\times d\hat{\tilde{R}}(x, z) \, d\tilde{F}(x_0, z_0)$$

(5)
$$+ o\{a^{(-1-\delta)/2}(\log a)^{m/2}\} \int (1+x_0) d\tilde{F}(x_0, z_0)$$

$$= \iint_{\substack{a < x \leq a+x_0 \\ x_0 < a/2 \\ |z_0| < \sqrt{a}}} [f(x-a, z) + O(1)\omega_f(x-a, z; e^{-\eta a})]$$
$$\times \hat{\tilde{r}}(x-x_0, z-z_0) \, d\tilde{\lambda}(x, z) \, d\tilde{F}(x_0, z_0)$$

$$+ o\{a^{(-1-\delta)/2}(\log a)^{m/2}\}$$

as $a \to \infty$. In the continuous case, $\tilde{\lambda}$ is a multiple of Lebesgue measure, and this expression can be simplified using dominated convergence. The mixed case is similar but more troublesome and will be treated carefully. Details for the continuous case are omitted. To avoid problems associated with the counting measure factor in $\lambda$ we will use the following rounding argument. Let $\langle w \rangle$ be $w$ with the last coordinate rounded to the nearest integer. Then

(6)
$$\int h \, d\tilde{\lambda} = \int h(x, \Sigma^{-1/2}(w - \gamma x)) \, d\lambda(x, w)$$
$$= \int h(x, \Sigma^{-1/2}(\langle w \rangle - \gamma x)) \, dx \, dw$$
$$= (a/\nu)^{m/2} \sqrt{|\Sigma|} \int h(x+a, z_a \sqrt{a/\nu}) \, dx \, dz$$

where

$$z_a = z_a(x) = \sqrt{\nu/a} \Sigma^{-1/2} (\langle \gamma(a+x) + \sqrt{a/\nu} \Sigma^{1/2} z \rangle - \gamma(a+x)).$$

Note that $|z_a - z|$ is bounded by some multiple of $1/\sqrt{a}$, uniformly in $x$. Using (6), the integral in (5) is

(7)
$$\iint_{\substack{a < x \leq a+x_0 \\ x_0 < a/2 \\ |z_0| < \sqrt{a}}} [f(x, z_a \sqrt{a/\nu}) + O(1)\omega_f(x, z_a \sqrt{a/\nu})] \Lambda \, dx \, dz \, d\tilde{F}(x_0, z_0),$$

where
$$\Lambda = \frac{1}{\nu}\left(1 + \frac{x-x_0}{a}\right)^{m/2} \phi\left(\frac{z_a}{\sqrt{1+(x-x_0)/a}} - \frac{z_0}{\sqrt{a}}\sqrt{\frac{\nu}{1+(x-x_0)/a}}\right)$$
$$\times \left[1 + \frac{1}{\sqrt{a}}\sqrt{\frac{\nu}{1+(x-x_0)/a}}\mathbf{H}\left(\frac{z_a}{\sqrt{1+(x-x_0)/a}} - \frac{z_0}{\sqrt{a}}\sqrt{\frac{\nu}{1+(x-x_0)/a}}\right)\right].$$

Viewing $\Lambda$ as a function of $a^{-1/2}$ (with $z_a$ treated as a separate independent variable), by Taylor's theorem with Lagrange's form for the remainder,

(8) $$\Lambda = \Lambda_0 + \theta/(2a),$$

where $\theta$ is $\partial^2\Lambda/\partial^2(a^{-1/2})$ evaluated at an intermediate point between 0 and $a^{-1/2}$, and
$$\Lambda_0 = \frac{1}{\nu}\phi(z_a)\left[1 + \sqrt{\frac{\nu}{a}}\left(\mathbf{H}(z_a) + z_0 \cdot z_a\right)\right].$$

A little calculus shows that $|\theta|$ is bounded by a multiple of

(9) $$\{1 + |x-x_0| + z_0^2\}(1+|z_a|^5)\phi(z_a)\exp\{|z_a|\sqrt{2\nu}\},$$

or by a multiple of this expression with $z_a$ changed to $z$. By (8),
$$a^{(1+\delta)/2}(\Lambda - \Lambda_0) \to 0$$

pointwise in $(x, z, x_0, z_0)$ as $a \to \infty$, and by (9), $a^{(1+\delta)/2}|\Lambda - \Lambda_0|$ is bounded by a multiple of
$$\{1 + |x-x_0|^{(1+\delta)/2} + |z_0|^{1+\delta}\}(1+|z|^5)\phi(z)\exp\{|z|\sqrt{2\nu}\}.$$

This last expression is integrable because
$$\int_{0<x\leq x_0}\left[1 + |x-x_0|^{(1+\delta)/2} + |z_0|^{1+\delta}\right]dx\, d\tilde{F}(x_0, z_0)$$
$$= E\left[X + \frac{2}{3+\delta}X^{(3+\delta)/2} + X|Z|^{1+\delta}\right] < \infty.$$

Hence by dominated convergence, (7) equals
$$\iint_{\substack{a<x\leq a+x_0 \\ x_0<a/2 \\ |z_0|<\sqrt{a}}} \left[f(x, z_a\sqrt{a/\nu}) + O(1)\omega_f(x, z_a\sqrt{a/\nu}; e^{-\eta a})\right]\Lambda_0\, dx\, dz\, d\tilde{F}(x_0, z_0)$$
$$+ o\{a^{(-1-\delta)/2}\}.$$

At this stage, by (3) and (4), there is no harm removing the upper endpoints from the $x_0$ and $z_0$ integrations. Performing the integral against $\tilde{F}$ and using (6),
$$Ef(X_\tau - a, Z_\tau)$$
$$= \int_{x>a}\left[f(x-a, z) + O(1)\omega_f(x-a, a; e^{-\eta a})\right]$$
$$\times \left\{\rho_0(x-a) + \sqrt{\nu/a}\left[\mathbf{H}(z\sqrt{\nu/a})\rho_0(x-a) + z\cdot\rho_1(x-a)\sqrt{\nu/a}\right]\right\}$$
$$\times \frac{\phi(z\sqrt{\nu/a})}{\sqrt{|\Sigma|}(a/\nu)^{m/2}}\, d\tilde{\lambda}(x, z)$$
$$+ o\{a^{(-1-\delta)/2}(\log a)^{m/2}\}$$

as $a \to \infty$. This proves Theorem 3.1. □

## 4. Ladder variables and the proof of Theorem 1.1

By introducing ladder variables, Theorem 1.1 can be obtained directly from the results of Section 3. Let $T_0 = 0$, $T_1 = T$ and define later ladder times recursively by
$$T_{k+1} = \inf\{n > T_k : X_n > X_{T_k}\}$$
for $k > 0$. Also, let $\tilde{X}_k = X_{T_k}$, $\tilde{W}_k = W_{T_k}$ and $\tilde{W}_k^* = W_{T_k}^*$. Then $\{(\tilde{X}_k, \tilde{W}_k)\}_{k \geq 0}$ and $\{(\tilde{X}_k, \tilde{W}_k^*)\}_{k \geq 0}$ are random walks. Define
$$\tilde{\tau} = \tilde{\tau}_a = \inf\{k : \tilde{X}_k > a\}.$$
Since $\tau$ must be a ladder time, $X_\tau = \tilde{X}_{\tilde{\tau}}$, $W_\tau = \tilde{W}_{\tilde{\tau}}$ and $W_\tau^* = \tilde{W}_{\tilde{\tau}}^*$. Hence Theorem 1.1 follows immediately from Theorem 3.1. The only work necessary is to show that $(\tilde{X}, \tilde{W})$ and $(\tilde{X}, \tilde{W}^*)$ satisfy the regularity conditions for Theorem 3.1, and to verify that the approximate densities in the two theorems agree. The necessary identities appear in Corollary 4.3 and Lemma 4.4.

The first result of this section concerns moments of randomly stopped sums. Let $e_1, e_2, \ldots$ be i.i.d. with $Ee_i = 0$ and $Ee_i^2 = 1$. Also, let the $e_i$ be independently adapted to a filtration $\mathcal{F}$, i.e., for each $n \geq 1$, $e_n$ is $\mathcal{F}_n$ measurable and independent of $\mathcal{F}_{n-1}$. Let $\zeta_n = \sum_{j=1}^n e_j$ and let $t$ be a stopping time.

**Lemma 4.1.** *Let $p > 0$. Then for some constant $K$ (that depends only on $p$),*

(10)  $$E|\zeta_t|^p \leq K\{Et^{p/2} + EtE|e_1|^p\}.$$

*If $Et < \infty$ then*

(11)  $$E[\zeta_t - \zeta_n; t > n] = 0$$

*and*

(12)  $$E[(\zeta_t - \zeta_n)^2; t > n] = E[t - n; t > n].$$

*If $Et^{3/2} < \infty$ and $E|e_1|^3 < \infty$ then*

(13)  $$E\zeta_t^3 = 3Et\zeta_t + EtEe_1^3.$$

*Proof.* Equations (11) and (12) are Wald's first and second identities. Under more stringent moment conditions, (13) was discovered by Chow, Robbins, and Teicher [5]. The moment conditions are weakened in Theorem 2.1 of Brown [2], and (13) follows from (10) using this result. For a narrower class of stopping times, (10) is given by Gut [8]. He notes (at the end of his paper) that the result can be easily obtained from Theorem 21.1 of Burkholder [3] (or Theorem 5.3 of [4]), and that approach works in this setting as well. □

For the next lemma, let $(e_1^{(1)}, e_1^{(2)})$, $(e_2^{(1)}, e_2^{(2)})$, ... be i.i.d., zero mean and independently adapted to $\mathcal{F}$. Let $\zeta_n^{(1)} = \sum_{i=1}^n e_i^{(1)}$, let $\zeta_n^{(2)} = \sum_{i=1}^n e_i^{(2)}$ and let $t$ be a stopping time.

**Lemma 4.2.** *If $E|e_1^{(1)}|^{3/2} < \infty$, $E|e_1^{(2)}|^3 < \infty$ and $Et^{3/2} < \infty$, then*
$$E\zeta_t^{(1)}\zeta_t^{(2)} = Et\, Ee_1^{(1)}e_1^{(2)}.$$

*Proof.* Since $\{\zeta_n^{(1)}\zeta_n^{(2)} - nEe_1^{(1)}e_1^{(2)}\}_{n\geq 1}$ is a zero mean martingale, it is sufficient to show that

$$E[\zeta_t^{(1)}\zeta_t^{(2)} - tEe_1^{(1)}e_1^{(2)} - \zeta_n^{(1)}\zeta_n^{(2)} + nEe_1^{(1)}e_1^{(2)}; t > n] \to 0.$$

This expression is the sum of two terms:

$$E[n-t; t > n]Ee_1^{(1)}e_1^{(2)} \to 0$$

and

(14) $$E[\zeta_t^{(1)}\zeta_t^{(2)} - \zeta_n^{(1)}\zeta_n^{(2)}; t > n].$$

Using Lemma 4.1, it is easy to show that $E[\zeta_t^{(1)} - \zeta_n^{(1)}|\mathcal{F}_n]$ is zero on $t > n$. By Hölder's inequality and (10), $E|(\zeta_t^{(1)} - \zeta_n^{(1)})\zeta_n^{(2)}| < \infty$, and so, conditioning on $\mathcal{F}_n$, $E[(\zeta_t^{(1)} - \zeta_n^{(1)})\zeta_n^{(2)}; t > n] = 0$. Similarly, $E[(\zeta_t^{(2)} - \zeta_n^{(2)})\zeta_n^{(1)}; t > n] = 0$. Hence, (14) equals

$$E\big[(\zeta_t^{(1)} - \zeta_n^{(1)})(\zeta_t^{(2)} - \zeta_n^{(2)}); t > n\big],$$

which is bounded in magnitude by

$$\left\{E|\zeta_t^{(1)} - \zeta_n^{(1)}|^{3/2}\right\}^{2/3}\left\{E\big[|\zeta_t^{(2)} - \zeta_n^{(2)}|^3; t > n\big]\right\}^{1/3}.$$

The second factor in this last expression approaches zero by the argument for (10). Using a square function inequality (such as Theorem 3.2 of [3]), $E|\zeta_t^{(1)} - \zeta_n^{(1)}|^{3/2}$ is bounded by a multiple of

$$E\left[\sum_{i=n}^{t} e_i^{(1)2}\right]^{3/4} \leq E\sum_{i=1}^{t}|e_i^{(1)}|^{3/2} = Et\,E|e_1^{(1)}|^{3/2},$$

and the lemma follows. □

**Corollary 4.3.** *Suppose $\delta \in (0,1)$, $\nu > 0$, $E|X|^{(3+\delta)/2} < \infty$ and $E|Z|^{3+\delta} < \infty$. Then $E\tilde{X}^{(3+\delta)/2} < \infty$, $ET^{(3+\delta)/2} < \infty$ and $E|\tilde{Z}|^{3+\delta} < \infty$. Also, $E\tilde{X} = \nu ET$, $E\tilde{W}/E\tilde{X} = \gamma$, $\mathrm{Cov}(\tilde{W} - \gamma\tilde{X}) = ET\Sigma$,*

(15) $$E(q \cdot \tilde{Z})^3 = \frac{E(q \cdot Z)^3}{\sqrt{ET}} + \frac{3|q|^2 ETq \cdot \tilde{Z}}{ET},$$

(16) $$E\tilde{Z}^2 q \cdot \tilde{Z} = \frac{EZ^2 q \cdot Z}{\sqrt{ET}} + \frac{(m+2)ETq \cdot \tilde{Z}}{ET}$$

*and*

(17) $$E(\tilde{X}/\nu - T)\tilde{Z} = \sqrt{ET}EXZ/\nu.$$

*If the moment condition for $X$ is strengthened to $E|X|^{3+\delta} < \infty$, then $E\tilde{X}^{3+\delta} < \infty$, $ET^{3+\delta} < \infty$, $E|\tilde{Z}^*|^{3+\delta} < \infty$, $E\tilde{W}^*/E\tilde{X} = \gamma^*$, $\mathrm{Cov}(\tilde{W}^* - \gamma^*\tilde{X}) = ET\Sigma_*$ and identities (15) and (16) relating $Z$ to $\tilde{Z}$ hold relating $Z^*$ to $\tilde{Z}^*$ after changing $m$ to $m+1$.*

*Proof.* By Theorem 2.1 of Gut [8], if $p \geq 1$, $EX > 0$ and $E|X|^p < \infty$, then $ET^p < \infty$ and $E\tilde{X}^p < \infty$. So (17) follows from Lemma 4.2 with $\zeta_n^{(1)} = X_n/\nu - n$ and $\zeta_n^{(2)}$ an arbitrary coordinate of $Z_n$. The other assertions follow easily from Lemma 4.1. For instance, by (12) with $n = 0$,

$$\mathrm{Var}\big(q \cdot (\tilde{W} - \gamma \tilde{X})\big) = E(q \cdot (W_T - \gamma X_T))^2$$
$$= ET \, \mathrm{Var}\big(q \cdot (W - \gamma X)\big).$$

So $q \cdot \mathrm{Cov}(\tilde{W} - \gamma \tilde{X})q = ETq \cdot \Sigma q$, which implies $\mathrm{Cov}(\tilde{W} - \gamma \tilde{X}) = ET\Sigma$. To obtain (16), let $\frac{1}{6}\sum_i \partial^2/\partial q_i^2$ act on (15). □

**Lemma 4.4.** *Suppose $EX > 0$ and $(X, W)$ satisfies Cramér's condition. Then $T$ is arithmetic on $\mathbb{Z}$ and $(\tilde{X}, \tilde{W})$ is strongly nonlattice with $T$.*

*Proof.* Introduce $T_- = \inf\{n : X_n \leq 0\}$. In the notation of Greenwood and Shaked [7], $T$ and $T_-$ are *dual* stopping times, and by their multivariate extension of the Wiener-Hopf factorization theorem,

$$1 - Ee^{i(\xi_1 X + \xi_2 \cdot W + \xi_3)}$$
$$= \left\{1 - Ee^{i(\xi_1 \tilde{X} + \xi_2 \cdot \tilde{W} + \xi_3 T)}\right\}\left\{1 - E\left[e^{i(\xi_1 X_{T_-} + \xi_2 \cdot W_{T_-} + \xi_3 T_-)}; T_- < \infty\right]\right\}.$$

The lemma now follows easily since $P(T_- < \infty) < 1$. □

## 5. Marginal distributions

To organize the proof of Theorem 1.2, it is convenient to approximate several marginal distributions separately. Lemma 5.1 gives marginal distributions for $W_\tau$ and $W_\tau^*$. Corollary 5.2 gives the marginal distribution for the first component of $W_\tau$, which establishes Theorem 1.2 when $h$ is linear. Lemma 2.2 shows that variables sufficiently close in a natural sense have the same asymptotic expansion. Lemma 5.3 provides an expansion for a special class of variables close to the variables in Theorem 1.2.

Let $\lambda_0$ be Lebesgue measure on $\mathbb{R}^m$ and $\lambda_0^*$ be the product of $\lambda_0$ with counting measure on $\mathbb{Z}$. Define $\hat{Q}_0$ and $\hat{Q}_0^*$ by

$$\frac{d\hat{Q}_0}{d\lambda_0}(w) = \frac{\phi(\tilde{q})}{\sqrt{|\Sigma|}(a/\nu)^{m/2}}\big\{1 + \sqrt{\nu/a}\mathbf{H}_0(\tilde{q})\big\}$$

and

$$\frac{d\hat{Q}_0^*}{d\lambda_0^*}(w^*) = \frac{\phi(\tilde{q}^*)}{\sqrt{|\Sigma_*|}(a/\nu)^{(m+1)/2}}\big\{1 + \sqrt{\nu/a}\mathbf{H}_0^*(\tilde{q}^*)\big\},$$

where

$$\tilde{q} = \Sigma^{-1/2}(w - \gamma a)\sqrt{\nu/a},$$
$$\tilde{q}^* = \Sigma_*^{-1/2}(w^* - \gamma^* a)\sqrt{\nu/a},$$
$$\mathbf{H}_0(\tilde{q}) = \frac{1}{6}E(\tilde{q} \cdot Z)^3 - \frac{1}{2}EZ^2 \tilde{q} \cdot Z$$
$$+ \frac{m + 2 - |\tilde{q}|^2}{2\nu}EX\tilde{q} \cdot Z + \frac{\tilde{q} \cdot \Sigma^{-1/2}\gamma}{2\nu ET}E\tilde{X}^2$$

and

$$\mathbf{H}_0^*(\tilde{q}^*) = \frac{1}{6}E(\tilde{q}^* \cdot Z^*)^3 - \frac{1}{2}EZ^{*2}\tilde{q}^* \cdot Z^*$$
$$+ \frac{m+3-|\tilde{q}^*|^2}{2\nu}EX\tilde{q}^* \cdot Z^* + \frac{\tilde{q}^* \cdot \Sigma_*^{-1/2}\gamma^*}{2\nu ET}E\tilde{X}^2.$$

**Lemma 5.1.** *Suppose $(X,W)$ satisfies Cramér's condition, $\nu > 0$, $E|X|^{2+\delta} < \infty$ and $E|Z|^{3+\delta} < \infty$, where $\delta \in (0,1)$. Then for some $\eta > 0$,*

$$(18) \quad Ef(W_\tau) = \int f \, d\hat{Q}_0 + O(1)\int \omega_f(\cdot; e^{-\eta a}) \, d\hat{Q}_0 + o\{a^{(-1-\delta)/2}(\log a)^{m/2}\}$$

*as $a \to \infty$, uniformly for nonnegative measurable $f$ bounded above by one. If the moment condition for $X$ is strengthened to $E|X|^{3+\delta}$, then for some $\eta > 0$,*

$$(19) \quad Ef(W_\tau^*) = \int f \, d\hat{Q}_0^* + O(1)\int \omega_f(\cdot; e^{-\eta a}) \, d\hat{Q}_0^* + o\{a^{(-1-\delta)/2}(\log a)^{(m+1)/2}\}$$

*as $a \to \infty$, uniformly for nonnegative measurable $f$ bounded above by one.*

*Proof.* Using Theorem 1.1, the only difficult task verifying (18) is to show that

$$(20) \quad \int f(w) \, d\hat{Q}(x,w) = \int f(w) \, d\hat{Q}_0(w) + o\{a^{(-1-\delta)/2}(\log a)^{m/2}\}.$$

Integrating over $w$, the $\hat{Q}$ marginal density for $x$ is bounded by some multiple of $\rho_0(x) + |\rho_1(x)|/\sqrt{a}$. Since

$$\int_{\sqrt{a}}^\infty \rho_0(x) \, dx = \frac{1}{\nu ET}E(\tilde{X} - \sqrt{a})^+$$
$$\leq \frac{a^{(-1-\delta)/2}}{\nu ET}E[\tilde{X}^{2+\delta}; \tilde{X} > \sqrt{a}]$$
$$= o\{a^{(-1-\delta)/2}\}$$

and

$$\frac{1}{\sqrt{a}}\int_{\sqrt{a}}^\infty |\rho_1(x)| \, dx \leq \frac{1}{\nu\sqrt{aET}}E[|\tilde{Z}|\tilde{X}; \tilde{X} > \sqrt{a}]$$
$$\leq \frac{a^{(-1-\delta)/2}}{\nu\sqrt{ET}}E[|\tilde{Z}|\tilde{X}^{1+\delta}; \tilde{X} > \sqrt{a}]$$
$$= o\{a^{(-1-\delta)/2}\},$$

$\hat{Q}\{[\sqrt{a},\infty) \times \mathbb{R}^m\} = o\{a^{(-1-\delta)/2}\}$. Hence there is no harm restricting the domain of integration in (20) to $x < \sqrt{a}$. Now $q - \tilde{q} = -\Sigma^{-1/2}\gamma x\sqrt{\nu/a}$, so any intermediate point on the line segment from $q$ to $\tilde{q}$ equals

$$\tilde{q} - \theta\Sigma^{-1/2}\gamma x\sqrt{\nu/a}$$

for some $\theta \in [0,1]$. The squared length of this vector is

$$|\tilde{q}|^2 - 2\theta\tilde{q} \cdot \Sigma^{-1/2}\gamma x\sqrt{\nu/a} + \theta^2 x^2 \gamma \cdot \Sigma^{-1}\gamma\nu/a \geq |\tilde{q}|^2 - 2|\tilde{q} \cdot \Sigma^{-1/2}\gamma x\sqrt{\nu/a}|.$$

Using Taylor's theorem with Lagrange's form for the remainder,

$$|\phi(q) - \phi(\tilde{q})| \leq \left\{|\tilde{q} \cdot \Sigma^{-1/2}\gamma\sqrt{\nu}| + \gamma \cdot \Sigma^{-1}\gamma x\nu/\sqrt{a}\right\} \frac{x}{\sqrt{a}} \phi(\tilde{q}) \exp|\tilde{q} \cdot \Sigma^{-1/2}\gamma x\sqrt{\nu/a}|$$

and

$$|\phi(q) - \phi(\tilde{q}) - \tilde{q} \cdot \Sigma^{-1/2}\gamma x\sqrt{\nu/a}\phi(\tilde{q})|$$
$$\leq \frac{1}{2}\left\{\left[|\tilde{q} \cdot \Sigma^{-1/2}\gamma\sqrt{\nu}| + \gamma \cdot \Sigma^{-1}\gamma x\nu/\sqrt{a}\right]^2 + \gamma \cdot \Sigma^{-1}\gamma\nu\right\}$$
$$\times \frac{x^2}{a}\phi(\tilde{q}) \exp|\tilde{q} \cdot \Sigma^{-1/2}\gamma x\sqrt{\nu/a}|.$$

Also,

$$|q \cdot \rho_1(x) - \tilde{q} \cdot \rho_1(x)| \leq |\rho_1(x)|\,|\Sigma^{-1/2}\gamma|x\sqrt{\nu/a},$$

and

$$|\mathbf{H}(\tilde{q}) - \mathbf{H}(q)| \leq K\big(|\tilde{q}|^3 + 1\big)x/\sqrt{a},$$

for some $K > 0$ on $x < \sqrt{a}$. Using these bounds, for an appropriate constant $K$,

$$\frac{d\hat{Q}}{d\lambda}(x,w)$$
$$= \frac{\phi(\tilde{q})}{\sqrt{|\Sigma|}(a/\nu)^{m/2}}\left\{\rho_0(x) + \sqrt{\nu/a}\big[\mathbf{H}(\tilde{q})\rho_0(x) + \tilde{q} \cdot \rho_1(x) + \tilde{q} \cdot \Sigma^{-1/2}\gamma x\rho_0(x)\big]\right\}$$
$$+ (1 + |\tilde{q}|^4)\big(x^{1+\delta}\rho_0(x) + x^\delta|\rho_1(x)|\big)\frac{\phi(\tilde{q})e^{K|\tilde{q}|}}{a^{m/2}}o\big\{a^{(-1-\delta)/2}\big\}$$

as $a \to \infty$, pointwise in $\tilde{q}$ and $x$. With $o$ changed to $O$, this result holds as $a \to \infty$, uniformly for $w \in \mathbb{R}^m$ and $x \leq \sqrt{a}$. Hence, by dominated convergence

$$\int f\,d\hat{Q} = \int \left\{\rho_0(x) + \sqrt{\nu/a}\big[\mathbf{H}(\tilde{q})\rho_0(x) + \tilde{q} \cdot \rho_1(x) + \tilde{q} \cdot \Sigma^{-1/2}\gamma x\rho_0(x)\big]\right\}$$
$$\times \frac{f(w)\phi(\tilde{q})}{\sqrt{|\Sigma|}(a/\nu)^{m/2}}\,d\lambda(x,w)$$
$$+ o\big\{a^{(-1-\delta)/2}\big\}$$
$$= \int \left\{1 + \sqrt{\nu/a}\bigg[\mathbf{H}(\tilde{q}) + \frac{E\tilde{X}\tilde{q}\cdot\tilde{Z}}{\nu\sqrt{ET}} + \frac{\tilde{q}\cdot\Sigma^{-1/2}\gamma}{2\nu ET}E\tilde{X}^2\bigg]\right\}$$
$$\times \frac{f(w)\phi(\tilde{q})}{\sqrt{|\Sigma|}(a/\nu)^{m/2}}\,d\lambda_0(w)$$
$$+ o\big\{a^{(-1-\delta)/2}\big\}$$

as $a \to \infty$. This gives (18). The proof of (19) is similar. □

When $m = 1$, integration of $d\hat{Q}_0/d\lambda_0$ gives the following corollary.

**Corollary 5.2.** *If $m = 1$, $(X, W)$ satisfies Cramér's condition, $\nu > 0$, $E|X|^{2+\delta} < \infty$ and $|Z|^{3+\delta} < \infty$, where $\delta \in (0,1)$, then*

$$P(W_\tau < w) = \Phi(\hat{w}) + \sqrt{\nu/a}\phi(\hat{w})\mathbf{H}_1(\hat{w}) + o\big\{a^{(-1-\delta)/2}\sqrt{\log a}\big\}$$

*as $a \to \infty$, uniformly in $w$, where $\hat{w} = (w - \gamma a)/(\sigma\sqrt{a/\nu})$, $\sigma^2 = \Sigma = \mathrm{Var}\,(W - \gamma X)$, and*

$$\mathbf{H}_1(\hat{w}) = (\hat{w}^2 - 1)\left\{-\frac{1}{6}EZ^3 + \frac{EXZ}{2\nu}\right\} - \frac{\gamma E\tilde{X}^2}{2\nu\sigma ET}.$$

For the next result, let $Y_n$ denote the first coordinate of $W_n$, so $W_n = (Y_n, V_n)$ where $V_n \in \mathbb{R}^{m-1}$. Also, let $V_n^* = (V_n, n)$, so $W_n^* = (Y_n, V_n, n) = (Y_n, V_n^*)$. Partition $\Sigma_*$ and $\gamma^*$ as

$$\Sigma_* = \begin{pmatrix} \Sigma_{11} & \Sigma_{12} \\ \Sigma_{21} & \Sigma_{22} \end{pmatrix} \quad \text{and} \quad \gamma^* = \begin{pmatrix} \gamma_1 \\ \gamma_2 \end{pmatrix},$$

where $\Sigma_{11}$ is 1 by 1, $\Sigma_{22}$ is $m$ by $m$, $\gamma_1 \in \mathbb{R}$ and $\gamma_2 \in \mathbb{R}^m$. Let

$$\Sigma_{11\cdot 2} = \Sigma_{11} - \Sigma_{12}\Sigma_{22}^{-1}\Sigma_{21}$$

and

$$\Sigma_{22\cdot 1} = \Sigma_{22} - \Sigma_{21}\Sigma_{11}^{-1}\Sigma_{12}.$$

The two factorizations of the multivariate normal density as a marginal times a conditional density give

$$\text{(21)} \quad \frac{\phi(\tilde{q}^*)}{\sqrt{|\Sigma|}(a/\nu)^{(m+1)/2}} = \frac{\phi(\tilde{v}^*)}{\sqrt{|\Sigma_{22}|}(a/\nu)^{m/2}} \frac{\phi(\tilde{y}_0 - r_1\tilde{v}^*)}{\sqrt{\Sigma_{11\cdot 2}a/\nu}}$$

$$= \frac{\phi(\tilde{y})}{\sqrt{\Sigma_{11}a/\nu}} \frac{\phi(\tilde{v}_0 - r_2\tilde{y})}{\sqrt{|\Sigma_{22\cdot 1}|}(a/\nu)^{m/2}}$$

where $w^* = (y, v)$,

$$r_1 = \Sigma_{11\cdot 2}^{-1/2}\Sigma_{12}\Sigma_{22}^{-1/2}, \qquad r_2 = \Sigma_{22\cdot 1}^{-1/2}\Sigma_{21}\Sigma_{11}^{-1/2}$$

$$\tilde{y}_0 = \Sigma_{11\cdot 2}^{-1/2}(y - \gamma_1 a)\sqrt{\nu/a}, \qquad \tilde{v}_0^* = \Sigma_{22\cdot 1}^{-1/2}(v^* - \gamma_2 a)\sqrt{\nu/a}$$

and

$$\tilde{y} = \Sigma_{11}^{-1/2}(y - \gamma_1 a)\sqrt{\nu/a}, \qquad \tilde{v}^* = \Sigma_{22}^{-1/2}(v^* - \gamma_2 a)\sqrt{\nu/a}.$$

Let $Y^\# = (Y_\tau - \gamma_1 a)/\sqrt{a}$ and $V^\# = (V_\tau^* - \gamma_2 a)/\sqrt{a}$. Finally, let $\lambda_1$ denote the product of Lebesgue measure on $\mathbb{R}^{m-1}$ with counting measure on $\mathbb{Z}$.

**Lemma 5.3.** *Let $h_0 \in \mathbb{R}$, $h_1 : \mathbb{R}^m \to \mathbb{R}$ a homogeneous linear function and $h_2 : \mathbb{R}^m \to \mathbb{R}$ a homogeneous quadratic function given by $h_2(q) = q \cdot Aq$ for some symmetric $m$ by $m$ matrix $A$. Assume $(X, W)$ satisfies Cramér's condition, $\nu > 0$, $E|X|^{3+\delta} < \infty$ and $E|Z|^{3+\delta} < \infty$, where $\delta \in (0, 1)$. Define*

$$\Xi_0 = \frac{Y^\# + [Y^\# h_1(V^\#) + h_2(V^\#)]/\sqrt{a}}{1 - h_0 Y^\#/\sqrt{a}}.$$

*Then*

$$P(\Xi_0 \leq c) = F_a(c) + o\{a^{(-1-\delta)/2}(\log a)^{(m+1)/2}\}$$

*as $a \to \infty$ uniformly in $c$, where $F_a$ is the approximate distribution function in Theorem 1.2.*

*Proof.* By Theorem 1.1, for any $\epsilon > 0$,

$$P(|Y^\#| \geq \epsilon\sqrt{a}) = o\{a^{(-1-\delta)/2}(\log a)^{(m+1)/2}\}$$

and

$$P(|V^\#| \geq \epsilon\sqrt{a}) = o\{a^{(-1-\delta)/2}(\log a)^{(m+1)/2}\}$$

as $a \to \infty$. Hence

$$P(\Xi_0 \leq c) = P\left(Y^\# \leq \frac{c - h_2(V^\#)/\sqrt{a}}{1 + [ch_0 + h_1(V^\#)]/\sqrt{a}}\right) + o\{a^{(-1-\delta)/2}(\log a)^{(m+1)/2}\}$$

as $a \to \infty$, uniformly if $|c|/\sqrt{a}$ stays sufficiently small. Using Lemma 5.1, this probability to the order of accuracy desired can be obtained by integrating the density of $\hat{Q}_0^*$ over the appropriate set. Using the first factorization of (21), integration over $y$ gives

$$
\begin{aligned}
(22) \quad \int & \frac{\phi(\tilde{v}^*)}{\sqrt{|\Sigma_{22}|}(a/\nu)^{m/2}} \\
& \times \left\{ \Phi(c_{v^*}) + \sqrt{\nu/a}\big[p_0(\tilde{v}^*)\Phi(c_{v^*}) + p_1(\tilde{v}^*, c_{v^*})\phi(c_{v^*})\big] \right\} d\lambda_1(v^*),
\end{aligned}
$$

where $p_0$ and $p_1$ are polynomials (their exact form is not important) and

$$
c_{v^*} = \sqrt{\frac{\nu}{\Sigma_{11\cdot 2}}}\, \frac{c - h_2(\Sigma_{22}^{1/2}\tilde{v}^*/\sqrt{\nu})/\sqrt{a}}{1 + \big[ch_0 + h_1(\Sigma_{22}^{1/2}\tilde{v}^*/\sqrt{\nu})\big]/\sqrt{a}} - r_1 \tilde{v}^*.
$$

As $a \to \infty$,

$$
c_{v^*} = \hat{c}_{v^*} - \chi(\tilde{v}^*)/\sqrt{a} + O(1/a)
$$

pointwise in $\tilde{v}^*$ and $\hat{c}$, where

$$
\hat{c}_{v^*} = c\sqrt{\nu/\Sigma_{11\cdot 2}} - r_1 \tilde{v}^*
$$

and

$$
\chi(\tilde{v}^*) = \sqrt{\nu/\Sigma_{11\cdot 2}}\Big\{ c^2 h_0 + c h_1(\Sigma_{22}^{1/2}\tilde{v}^*/\sqrt{\nu}) + h_2(\Sigma_{22}^{1/2}\tilde{v}^*/\sqrt{\nu}) \Big\}.
$$

Then

$$
\Phi(c_{v^*}) = \Phi(\hat{c}_{v^*}) - \phi(\hat{c}_{v^*})\chi(\tilde{v}^*)/\sqrt{a} + O(1/a),
$$

$$
\Phi(c_{v^*})/\sqrt{a} = \Phi(\hat{c}_{v^*})/\sqrt{a} + O(1/a)
$$

and

$$
p_1(\tilde{v}^*, c_{v^*})\phi(c_{v^*})/\sqrt{a} = p_1(\tilde{v}^*, \hat{c}_{v^*})\phi(\hat{c}_{v^*})/\sqrt{a} + O(1/a)
$$

as $a \to \infty$, pointwise in $\tilde{v}^*$ and $\hat{c}$. For some $\epsilon > 0$, with the error rate changed to $o\{a^{(-1-\delta)/2}\}$, these three relations will hold uniformly for $|\tilde{v}^*| < a^\epsilon$ and $|\hat{c}| < a^\epsilon$. This is sufficient to justify the obvious substitutions into (22). This has a net effect of changing $c_{v^*}$ in (22) to $\hat{c}_{v^*}$ and adding an additional term:

$$
(23) \quad -\int \frac{\phi(\tilde{v}^*)\phi(\hat{c}_{v^*})}{\sqrt{|\Sigma_{22}|}(a/\nu)^{m/2}} \frac{\chi(\tilde{v}^*)}{\sqrt{a}}\, d\lambda_1(v^*).
$$

Without this extra term, the integral is simply the approximation for $P(W^\# \le c)$. To finish it is sufficient to show that (23) equals

$$
-\phi(\hat{c})\sqrt{\frac{\nu}{a}}\left\{ \frac{\sigma \hat{c}^2 h_0}{\nu} + \frac{\hat{c}^2 h_1(\Sigma_{21})}{\nu \sigma} + \frac{\hat{c}^2 - 1}{\nu \sigma^3}\Sigma_{12} A \Sigma_{21} + \frac{1}{\nu \sigma}\operatorname{tr}\big[A\Sigma_{22}\big] \right\} + o\{a^{(-1-\delta)/2}\},
$$

as this is the difference between $F_a(c)$ and the approximation for $P(W^\# \le c)$ obtained from Corollary 5.2. Using both factorizations in (21),

$$
\frac{\phi(\tilde{v}^*)\phi(\hat{c}_{v^*})}{\sqrt{|\Sigma_{22}|}(a/\nu)^{m/2}} = \frac{\phi(\hat{c})}{\sqrt{\Sigma_{11} a/\nu}} \frac{\phi(\tilde{v}_0^* - r_2 \hat{c})}{\sqrt{|\Sigma_{22\cdot 1}|}(a/\nu)^{m/2}} \sqrt{\Sigma_{11\cdot 2} a/\nu}.
$$

In the integral in (23) there is no harm changing the measure of integration from $\lambda_1$ to Lebesgue measure. This can be verified using the Euler–MacLaurin sum formula: If $|f'|$ is integrable, then

$$\sum_{n \in \mathbb{Z}} f(n) = \int f(x)\, dx + \int \langle x \rangle_0 f'(x)\, dx$$

where $\langle \cdot \rangle_0$ is the periodic function defined for $x \in [0,1)$ by $\langle x \rangle_0 = x - 1/2$. In (23), partial derivatives with respect to the last component of $v^*$ introduce an extra factor of $\sqrt{a}$ in the denominator, so changing $\lambda_1$ to Lebesgue measure changes the integral by $O(1/a)$. Since $\phi(\tilde{v}_0^* - r_2 \hat{c})/(\sqrt{|\Sigma_{22 \cdot 1}|}(a/\nu)^{m/2})$ is a multivariate normal density, (23) has been simplified to the expectation of a quadratic function of a normal vector. The lemma then follows from standard facts about the normal distribution. Uniformity for all $c$ follows from uniformity for $|\hat{c}| \leq a^\epsilon$ by monotonicity of the distribution function being approximated. □

*Proof of Theorem 1.2.* Using Lemma 2.2, this proof will be accomplished showing that

$$\frac{a^{(1+\delta)/2}}{\sqrt{\log a}}(\Xi - \Xi_0) \to 0$$

in probability at rate $o\{a^{(-1-\delta)/2}\sqrt{\log a}\}$. Using Theorem 4 of Lorden [14], for some $K > 0$,

$$(24) \quad \begin{aligned} P\big(X_\tau - a \geq \epsilon a^{(1-\delta)/2}\big) &\leq K E\big[X; X \geq \epsilon a^{(1-\delta)/2}\big] \\ &= o\{a^{(-1-\delta)/2}\} \end{aligned}$$

as $a \to \infty$, since the conditions of Theorem 1.2 ensure $E|X|^{2/(1-\delta)} < \infty$. Hence

$$(25) \quad a^{(1-\delta)/2}(X_\tau - a) \to 0$$

in probability at rate $o\{a^{(-1-\delta)/2}\}$. Applying Theorem 1.1 with the random walk changed to $\{(X_n, Y_n)\}$, for some constant $K > 0$,

$$P\big(|Y^\#| \geq K\sqrt{\log a}\big) = o\{a^{(-1-\delta)/2}\sqrt{\log a}\}.$$

Similarly, applying Theorem 1.1 to two dimensional random walks where the second coordinates are projections of $W_n$, for some $K > 0$,

$$(26) \quad P\big(|W^\#| \geq K\sqrt{\log a}\big) = o\{a^{(-1-\delta)/2}\sqrt{\log a}\}.$$

Let

$$\Xi_1 = Y^\# + \big[h_0 Y^{\#^2} + Y^\# h_1(V^\#) + h_2(V^\#)\big]/\sqrt{a}.$$

By Taylor expansion, on $\{|W^\#| < K\sqrt{\log a}\}$,

$$\Xi_0 - \Xi_1 = O\{(\log a)^{3/2}/a\}.$$

Consequently,

$$a^{(1+\delta)/2}(\Xi_0 - \Xi_1) \to 0$$

in probability at rate $o\{a^{(-1-\delta)/2}\sqrt{\log a}\}$. Next, let

$$\Xi_2 = \Xi_1 + h_3(X_\tau/a - 1)^2 + (X_\tau/a - 1)h_4(Y_\tau/a - \gamma_1, V_\tau^*/a - \gamma_2),$$

the 2-term Taylor approximation for $\Xi$. Using (25) and (26),

$$(\Xi_2 - \Xi_1)a^{(1+\delta)/2}/\sqrt{\log a} \to 0$$

in probability at rate $o\{a^{(-1-\delta)/2}\sqrt{\log a}\}$. By (24) and (26),

$$P(S_\tau^*/a \notin N_0) = o\{a^{(-1-\delta)/2}\sqrt{\log a}\}$$

as $a \to \infty$, for any neighborhood $N_0$ of $s_0$. With $N_0$ sufficiently small and $K$ large enough, on $\{S_\tau^*/a \in N_0\}$

$$\Xi - \Xi_2 \leq K\sqrt{a}\{|X_\tau/a - 1|^3 + |W_\tau^*/a - \gamma|^3\}$$

by Taylor's theorem. Using (25) and (26), it follows that

$$(\Xi - \Xi_2)a^{(1+\delta)/2}/\sqrt{\log a} \to 0$$

in probability at rate $o\{a^{(-1-\delta)/2}\sqrt{\log a}\}$, which proves Theorem 1.2. □

## Acknowledgments

Portions of this research were accomplished while I was on sabbatical leave at Sydney University. Suggestions from a referee were useful and quite appreciated.